\documentclass[10pt]{article}
\usepackage{amssymb}
\usepackage{amsmath}
\usepackage[dvips]{epsfig}
\usepackage{array}
\newcommand{\Var}{\operatorname{Var}}
\newcommand{\SZ}{\operatorname{SZ}}
\newcommand{\DS}{\operatorname{DS}}
\newcommand{\WS}{\operatorname{WS}}
\newcommand{\DR}{\operatorname{DR}}
\newcommand{\WR}{\operatorname{WR}}
\newcommand{\DZ}{\operatorname{DZ}}
\newcommand{\WZ}{\operatorname{WZ}}
\newcommand{\XY}{\operatorname{XY}}
\newcommand{\XX}{\operatorname{XX}}
\newcommand{\SR}{\operatorname{SR}}
\newcommand{\ZR}{\operatorname{ZR}}
\renewcommand{\SS}{\operatorname{SS}}
\newcommand{\ZZ}{\operatorname{ZZ}}
\newcommand{\RR}{\operatorname{RR}}

\title{Calculating the correlation coefficients of graph-theoretical indices}
\author{Stephan G. Wagner\thanks{Department of Mathematics, Graz University of Technology, Steyrergasse 30, A-8010 Graz, Austria ({\tt wagner@finanz.math.tugraz.at}).}}
\date{}

\begin{document}

\maketitle

\begin{abstract}
Using a generating function approach, the correlation coefficients of four
different graph-theoretical indices, namely the number of independent vertex
subsets, the number of matchings, the number of subtrees and the Wiener index,
are asymptotically determined for random rooted ordered trees.
\end{abstract}

\pagestyle{plain}

\section{Introduction}

In \cite{wagcor}, the author investigates correlation measures for
graph-theoretical indices which are of interest in theoretical chemistry. In
particular, the correlation coefficients for these indices are asymptotically
determined. Since the necessary calculations are rather lengthy and tedious,
only a few of them are explicitly provided there. This additional note fills
the gap by explaining the involved details. See \cite{wagcor} for further
references and applications.

The underlying stochastic model is the following: of all rooted ordered trees
on $n$ vertices, a tree $T_n$ is selected uniformly at random. The parameter
we are interested in is the correlation coefficient of two indices $X_n =
X(T_n)$ and $Y_n = Y(T_n)$, defined by
\begin{equation}
r(X_n,Y_n) = \frac{E(X_nY_n)-E(X_n)E(Y_n)}{\sqrt{\Var(X_n)\Var(Y_n)}}.
\end{equation}
Here, $X$ and $Y$ are two of the following four indices:
\renewcommand{\labelenumi}{(\theenumi)}
\begin{enumerate}
\item The \emph{Merrifield-Simmons-} or $\sigma$-index is defined to be the
number of independent vertex subsets of a graph, i.e. the number of vertex
subsets in which no two vertices are adjacent, including the empty set.
Merrifield and Simmons investigated the $\sigma$-index in their work
\cite{mer89}
and pointed out its correlation to boiling points of molecules.
\item The \emph{Hosoya-} or $Z$-index (\cite{hos86}) is defined as the number
of independent edge subsets (also referred to as ``matchings''), i.e. the
number of edge subsets in which no two edges are adjacent, including the empty
set again.
\item The \emph{number of subtrees} is called $\rho$-index in \cite{mer89} and
was discussed lately in a paper of Sz\'ekely and Wang \cite{sze05}.
\item The \emph{Wiener index} is probably the most popular topological index
(s. \cite{dob01,ent94,wie47}). It is defined as the sum of all the distances
between pairs of vertices, i.e.
\begin{equation}
W(G) = \sum_{v,w \in V(G)} d_G(v,w).
\end{equation}
\end{enumerate}

It will be shown how one can derive functional equations for generating
functions yielding to asymptotic formulas for the quantities of interest by
means of the Flajolet-Odlyzko singularity analysis \cite{fla90}. It
turns out that the approach is a slightly different one for the correlation
with the Wiener index in view of its different growth structure. Therefore, the
correlation coefficients of the $\sigma$-, $Z$- and $\rho$-index are determined
in the following section, whereas the correlation to the Wiener index is
investigated in Section~\ref{wiener}.

\section{$\sigma$-, $Z$- and $\rho$-index}

In this section, we want to determine the asymptotic behavior of the
generating function
$$\sum_{T} X(T)Y(T)z^{|T|},$$
where $X,Y$ stand for $\sigma$-, $Z$- or $\rho$-index (possibly, $X$ and $Y$
are the same). They count, respectively, the number of independent vertex
subsets, edge subsets and subtrees of a tree $T$. This is done by
distinguishing between two cases for each of the indices:
\begin{itemize}
\item the root vertex belongs to the independent vertex subset/edge
  subset/subtree,
\item the root does not belong to it.
\end{itemize}
We denote the corresponding quantities by $\sigma_1(T),\sigma_2(T)$
resp. $Z_1(T),Z_2(T)$ and $\rho_1(T),\rho_2(T)$. If $T_1,\ldots,T_k$ are the
branches of the rooted tree $T$, it is easy to see that
\begin{align*}
&\sigma_1(T) = \prod_{i=1}^k \sigma_2(T_i), \\
&\sigma_2(T) = \prod_{i=1}^k (\sigma_1(T_i) + \sigma_2(T_i)), \\
\end{align*}
\begin{align*}
&Z_1(T) = \sum_{j=1}^k Z_2(T_j) \prod_{\underset{i \neq j}{i=1}}^k (Z_1(T_i) +
Z_2(T_i), \\
&Z_2(T) = \prod_{i=1}^k (Z_1(T_i) + Z_2(T_i)), \\
&\rho_1(T) = \prod_{i=1}^k (1+ \rho_1(T_i)), \\
&\rho_2(T) = \sum_{i=1}^k (\rho_1(T_i) + \rho_2(T_i)). \\
\end{align*}
The corresponding generating functions are called $S_1,S_2$ resp. $Z_1,Z_2$
and $R_1,R_2$. Functional equations for these functions which follow from the
recursive relations given above and lead to asymptotic formulas for the average
indices have already been presented by Klazar \cite{kla97} and others. For the
sake of completeness, these are given here as well:
\begin{align*}
S_1(z) &= \frac{z}{1-S_2(z)}, \\
S_2(z) &= \frac{z}{1-S_1(z)-S_2(z)}, \\
Z_1(z) &= \frac{zZ_2(z)}{(1-Z_1(z)-Z_2(z))^2}, \\
Z_2(z) &= \frac{z}{1-Z_1(z)-Z_2(z)}, \\
R_1(z) &= \frac{z}{1-R_1(z)-T(z)}, \\
R_2(z) &= \frac{z}{(1-T(z))^2}(R_1(z)+R_2(z)).
\end{align*}
Here, $T(z)$ is the generating function for the number of rooted ordered trees.
It is well-known that $T(z)$ satisfies the functional equation
$$T(z) = \frac{z}{1-T(z)},$$
which leads to an explicit formula for the number of rooted ordered trees:
$$t_n = \frac{1}{n} \binom{2n-2}{n-1} \sim \frac{1}{4\sqrt{\pi}} n^{-3/2} 4^n.$$
As determined by Klazar \cite{kla97}, the functional equations presented above
yield the following asymptotic formulas for the expected values of our indices:
$$E(\sigma_n) \sim \sqrt{3} \left( \frac{27}{16} \right)^{n-1} \approx
(1.02640) \cdot (1.6875)^n.$$
$$E(Z_n) \sim \sqrt{\frac{65-\sqrt{13}}{78}} \left( \frac{35+13\sqrt{13}}{54}
\right)^n \approx (0.88719) \cdot (1.51615)^n.$$
$$E(\rho_n) \sim \frac{16}{3\sqrt{15}} \left( \frac{25}{16} \right)^n \approx
(1.37706) \cdot (1.5625)^n.$$

If we combine $\sigma_{1,2}$ and $Z_{1,2}$, we
obtain generating functions $\SZ_{11},\SZ_{12},\SZ_{21},\SZ_{22}$ for the
correlation between $\sigma$- and $Z$-index:
$$\SZ_{ij}(z) = \sum_{T} \sigma_i(T)Z_j(T)z^{|T|}.$$
The total is denoted by $\SZ$. In an analogous manner, we define
$\SR_{ij}$,$\ZR_{ij}$,$\SS_{ij}$,$\ZZ_{ij}$ and $\RR_{ij}$.

Next, we observe that $1 \leq \sigma(T),Z(T),\rho(T) \leq 2^{|T|}$ for trivial
reasons. This helps us to restrict the range of the radius of convergence. In
all our cases, it must lie within the interval $\left[\frac{1}{16},\frac{1}{4}
\right]$. Even more can be told about it: the radius of
convergence of $\XY$ can be at most the minimum of the radii of $X$ and
$Y$ (which are $\frac{4}{27}$, $\frac{13\sqrt{13}-35}{72}$ and $\frac{4}{25}$
for $\sigma$-, $Z$- and $\rho$-index respectively). Furthermore, since
$\sigma(T) \geq F_{|T|+2}$, where $F_n$ denotes the $n$-th Fibonacci number (a
result due to Prodinger and Tichy \cite{pro82}), we even know that, for
example, the radius of convergence of $\SR$ is $\leq \frac{4}{25} \cdot
\frac{\sqrt{5}-1}{2} = \frac{2(\sqrt{5}-1)}{25}$. For the generating functions
$\SS$, $\ZZ$ and $\RR$, we may apply the Cauchy-Schwarz-inequality to obtain an
upper bound for the convergence radius easily: if $X$ has a convergence radius
$r$, then $\XX$ has convergence radius $\leq 4r^2$. Summing up, we have the
following bounds for the convergence radii:
\begin{itemize}
\item $\SZ$: interval
  $\left[\frac{1}{16},\frac{(13\sqrt{13}-35)(\sqrt{5}-1)}{144}\right]$,
\item $\SR$: interval
  $\left[\frac{1}{16},\frac{2(\sqrt{5}-1)}{25}\right]$,
\item $\ZR$: interval
  $\left[\frac{1}{16},\frac{4}{25}\right]$,
\item $\SS$: interval
  $\left[\frac{1}{16},\frac{64}{729}\right]$,
\item $\ZZ$: interval
  $\left[\frac{1}{16},\frac{1711-455\sqrt{13}}{648}\right]$,
\item $\RR$: interval
  $\left[\frac{1}{16},\frac{64}{625}\right]$.
\end{itemize}
In all the cases, we will see that these estimates are sufficient to
determine the correct dominating singularity.

\subsection{$\sigma$- and $Z$-index}

The recursive relations for $\sigma$- and $Z$-index lead to the following
system of functional equations:
\begin{equation*}
\begin{split}
\SZ_{11}(z) &= \frac{z\SZ_{22}(z)}{(1-\SZ_{21}(z)-\SZ_{22}(z))^2}, \\
\SZ_{12}(z) &= \frac{z}{1-\SZ_{21}(z)-\SZ_{22}(z)}, \\
\SZ_{21}(z) &= \frac{z(\SZ_{12}(z)+\SZ_{22}(z))}{(1-\SZ_{11}(z)-\SZ_{12}(z)-\SZ_{21}(z)-\SZ_{22}(z))^2}, \\
\SZ_{22}(z) &= \frac{z}{1-\SZ_{11}(z)-\SZ_{12}(z)-\SZ_{21}(z)-\SZ_{22}(z)}.
\end{split}
\end{equation*}
For instance, the functional equation for $\SZ_{11}$ can be deduced as follows:
\begin{align*}
\SZ_{11}(z) &= \sum_T \sigma_1(T) Z_1(T) z^{|T|} \\
&= \sum_{k \geq 0} \sum_{j=1}^k \sum_{T_1} \sum_{T_2} \ldots \sum_{T_k} \left(
\sigma_2(T_j) Z_2(T_j) \prod_{i \neq j} \sigma_2(T_i) (Z_1(T_i) + Z_2(T_i))
\right) \\
&\ \ \ \cdot z^{|T_1| + \ldots + |T_k| + 1} \\
&= z \sum_{k \geq 0} k \SZ_{22}(z) (\SZ_{21}(z) + \SZ_{22}(z))^{k-1} \\
&= \frac{z\SZ_{22}(z)}{(1-SZ_{21}(z) - \SZ_{22}(z))^2}.
\end{align*}
From these, a single equation for $\SZ_{22}$ can be worked out by
means of Gr\"obner bases \cite{fro97}, and this can be used to determine the
dominating singularity of $\SZ$. All computations are given in the
accompanying Mathematica\textsuperscript{\textregistered} files, which can be
found on
\texttt{http://finanz.math.tugraz.at/{\textasciitilde}wagner/Correlation}. From
the equation
$$s^{10} + 2z s^8 - 3z s^7 + z^2 s^6 - 4z^2 s^5
+ 3z^2 s^4 - z^3 s^3 + 2z^3 s^2 - z^3 s + z^4 = 0$$
that is satisfied by $s=\SZ_{22}(z)$, we find that
$$\SZ_{22}(z) \sim 0.171502 - 0.138532 \sqrt{1-\frac{z}{z_0}}$$
around the singularity $z_0 = 0.0982673$. From the relation $\SZ(z) = 1 -
\frac{z}{\SZ_{22}(z)}$ we obtain
$$\SZ(z) \sim 0.427020 - 0.462827 \sqrt{1-\frac{z}{z_0}},$$
which gives us the asymptotic formula
$$E(\sigma_nZ_n) \sim (0.92565) \cdot (2.54408)^n.$$
by a simple application of the Flajolet-Odlyzko singularity analysis
\cite{fla90}, as it is also explained in \cite{wagcor}.

\subsection{$\sigma$- and $\rho$-index}

The recursive relations for $\sigma$- and $\rho$-index lead to the following
system of functional equations:
\begin{equation*}
\begin{split}
\SR_{11}(z) &= \frac{z}{1-\SR_{21}(z)-S_2(z)}, \\
\SR_{12}(z) &= \frac{z(\SR_{21}(z) + \SR_{22}(z))}{(1-S_2(z))^2} =
\frac{S_1(z)^2(\SR_{21}(z) + \SR_{22}(z))}{z}, \\
\SR_{21}(z) &= \frac{z}{1-\SR_{11}(z)-\SR_{21}(z)-S_1(z)-S_2(z)}, \\
\SR_{22}(z) &= \frac{z(\SR_{11}(z) + \SR_{12}(z) + \SR_{21}(z) +
  \SR_{22}(z))}{(1-S_1(z)-S_2(z))^2} \\
&=  \frac{S_2(z)^2(\SR_{11}(z) + \SR_{12}(z) + \SR_{21}(z) + \SR_{22}(z))}{z}.
\end{split}
\end{equation*}
It would be possible to carry out the same procedure as in the previous
case; however, it saves computing time to consider $\SR_{11}$ and $\SR_{21}$
first. Then, $\SR_{12}$ and $\SR_{22}$ are given by simple linear equations
which result in the formula
$$\SR(z) = \frac{z(S_1(z)^2 \SR_{21}(z) + z\SR_{11}(z) +
  z\SR_{21}(z))}{z^2-zS_2(z)^2 - S_1(z)^2S_2(z)^2}.$$
We know that $S_1$ and $S_2$ are holomorphic in the region of interest, so the
dominating singularity of $\SR$ is either a singularity of $\SR_{11}$ and
$\SR_{21}$ or a zero of the denominator. However, using the functional
equations for $S_1$ and $S_2$, we find that the denominator only vanishes at $z
= 0$ and at $z = \frac{4}{27}$. Since $\frac{4}{27}$ does not lie in our
estimated interval, we have to determine the singularities of $\SR_{11}$ and
$\SR_{21}$. By the Gr\"obner basis approach, we find that $s = \SR_{11}(z)$
satisfies the polynomial equation
$$z s^9 - 6 z^2 s^7 - 4 z^3 s^6 + 7 z^3 s^5 - 2 z^4 s^4 - 
    (z^5 + 3 z^4) s^3 + z^5 s^2 + z^5 s - z^6 = 0.$$
0 is clearly not a singularity, and the other common zeroes of the polynomial
and its derivative satisfy the polynomial equation
$$(27z-4)(8z^2+81)^2(125z^3-412z^2-40936z+3844) = 0.$$
The only value that lies within our bounds and satisfies the equation is $z_0
\approx 0.0938166$. Since $\SR_{21}(z) = 1 - S_2(z) - \frac{z}{\SR_{11}(z)}$
and $\SR_{11}$ only vanishes at $z=0$, the smallest singularity of $\SR_{21}$
is the same. Expanding around $z_0$ gives us
$$\SR(z) \sim 0.560623 - 0.683264 \sqrt{1-\frac{z}{z_0}}$$
and thus
$$E(\sigma_n \rho_n) \sim (1.36653) \cdot (2.66477)^n.$$

\subsection{$Z$- and $\rho$-index}

The recursive relations for $Z$- and $\rho$-index lead to the following
system of functional equations:
\begin{equation*}
\begin{split}
\ZR_{11}(z) &= \frac{z(\ZR_{21}(z) + Z_2(z))}
{(1-\ZR_{11}(z)-\ZR_{21}(z)-Z_1(z)-Z_2(z))^2}, \\
\ZR_{12}(z) &= \frac{z(\ZR_{21}(z) + \ZR_{22}(z))}{(1-Z_1(z)-Z_2(z))^2} +
\frac{2zZ_2(z)\ZR(z)}{(1-Z_1(z)-Z_2(z))^3} \\
&= \frac{Z_2(z)^2(\ZR_{21}(z) +
  \ZR_{22}(z))}{z} + \frac{2Z_1(z)Z_2(z)\ZR(z)}{z}, \\
\ZR_{21}(z) &= \frac{z}{1-\ZR_{11}(z)-\ZR_{21}(z)-Z_1(z)-Z_2(z)}, \\
\ZR_{22}(z) &= \frac{z\ZR(z)}{(1-Z_1(z)-Z_2(z))^2} =  \frac{Z_2(z)^2\ZR(z)}{z}.
\end{split}
\end{equation*}
Using the same approach as in the previous example, we arrive at
$$\ZR(z) = \frac{z(Z_2(z)^2 \ZR_{21}(z) + z\ZR_{11}(z) +
  z\ZR_{21}(z))}{z^2-2zZ_1(z)Z_2(z) - zZ_2(z)^2-Z_2(z)^4}.$$
Again, $Z_1$ and $Z_2$ are holomorphic in the region of interest, and the
denominator only vanishes at $z = 0$ and at $z =
\frac{13\sqrt{13}-35}{72}$, which does not lie in our 
estimated interval, so we have to determine the singularities of $\ZR_{11}$ and
$\ZR_{21}$. We observe first that $\ZR_{11}(z) =
\frac{\ZR_{21}(z)^2(\ZR_{21}(z)+Z_2(z))}{z}$ and $Z_1(z) = \frac{Z_2(z)^3}{z}$. Inserting
yields
\begin{multline*}
z^2- z\ZR_{21}(z) + zZ_2(z)\ZR_{21}(z) + Z_2(z)^3\ZR_{21}(z) \\
+ z\ZR_{21}(z)^2 + Z_2(z)\ZR_{21}(z)^3 + \ZR_{21}(z)^4 = 0.
\end{multline*}
Elimination of $Z_2$ via the relation $z^2 - zZ_2(z) + zZ_2(z)^2 + Z_2(z)^4 =
0$ gives a polynomial equation of degree 16. For computational purposes,
however, it is much faster to find a common zero of the equation above together
with its derivative and the condition for $Z_2$.
Then we find that a singularity of $\ZR_{21}$ (and thus also of
$\ZR_{11}$) and $\ZR$ must satisfy the polynomial equation
$$(4096z^2-448z+1)(2560000z^2+2894400z+531441) = 0.$$
The only value that lies within our bounds and satisfies the equation is $z_0
\approx 0.107095$. Expanding around $z_0 = \frac{7+3\sqrt{5}}{128}$ gives us 
$$\ZR(z) \sim \frac{1}{928}(211+93\sqrt{5}) -
\frac{1}{232}\sqrt{\frac{5(128985+57683\sqrt{5})}{58}} \cdot
\sqrt{1-\frac{z}{z_0}}$$
and thus
$$E(Z_n \rho_n) \sim \frac{1}{116}\sqrt{\frac{5(128985+57683\sqrt{5})}{58}}
\cdot \left( 8(7-3\sqrt{5}) \right)^n.$$

\newpage

\subsection{Variance of the $\sigma$-index}

For the variances, the calculations are even a little simpler, since we have
one variable less to cope with. In this case, we obtain the functional
equations 
\begin{equation*}
\begin{split}
\SS_{11}(z) &= \frac{z}{1-\SS_{22}(z)}, \\
\SS_{12}(z) &= \frac{z}{1-\SS_{12}(z)-\SS_{22}(z)}, \\
\SS_{22}(z) &= \frac{z}{1-\SS_{11}(z)-2\SS_{12}(z)-\SS_{22}(z)},
\end{split}
\end{equation*}
which result in a single equation for $s = \SS(z) = \SS_{11}(z) + 2\SS_{12}(z) 
+ \SS_{22}(z)$:
\begin{multline*}
s^6-6s^5+(4z+14)s^4+(8z^2-20z-16)s^3+(4z^3-30z^2+36z+9)s^2 \\
-(12z^3-36z^2+28z+2)s-(z^4-8z^3+14z^2-8z) = 0.
\end{multline*}
We use Gr\"obner bases once again and find the only possible singularity
which lies within our bounds: its value is $z_0 \approx 0.0873832$. We
calculate the expansion of $\SS(z)$ around $z_0$:
$$\SS(z) \sim 0.614803 - 0.519010 \sqrt{1-\frac{z}{z_0}},$$
and finally obtain the asymptotics of the variance:
$$\Var(\sigma_n) \sim (1.03802) \cdot (2.86096)^n.$$

\subsection{Variance of the $Z$-index}

We proceed in the same way in the case of the $Z$-index:
\begin{equation*}
\begin{split}
\ZZ_{11}(z) &= \frac{z\ZZ_{22}(z)}{(1-\ZZ(z))^2} +
\frac{2z(\ZZ_{12}(z) + \ZZ_{22}(z))^2}{(1-\ZZ(z))^3}, \\
\ZZ_{12}(z) &= \frac{z(\ZZ_{12}(z) + \ZZ_{22}(z))}{(1-\ZZ(z))^2}, \\
\ZZ_{22}(z) &= \frac{z}{1-\ZZ(z)},
\end{split}
\end{equation*}
which results in a single equation for $s = \ZZ(z) = \ZZ_{11}(z) + 2\ZZ_{12}(z) 
+ \ZZ_{22}(z)$:
\begin{multline*}
s^8-7s^7-(z-21)s^6+(4z-35)s^5+(2z^2-5z+35)s^4 -(7z^2+21)s^3 \\
-(z^3-9z^2-5z-7)s^2+(2z^3-5z^2-4z-1)s+(z^4-z^3+z^2+z) = 0.
\end{multline*}
Here, the value of the singularity is $z_0 \approx 0.107969$. We
calculate the expansion of $\ZZ(z)$ around $z_0$:
$$\ZZ(z) \sim 0.296221 - 0.386136 \sqrt{1-\frac{z}{z_0}},$$
and finally obtain the asymptotics of the variance:
$$\Var(Z_n) \sim (0.77227) \cdot (2.31549)^n.$$

\subsection{Variance of the $\rho$-index}

In this case, we obtain the following system of equations:
\begin{equation*}
\begin{split}
\RR_{11}(z) &= \frac{z}{(1-\RR_{11}(z)-2R_1(z)-T(z)}, \\
\RR_{12}(z) &= \frac{z(\RR_{11}(z) +
  \RR_{12}(z)+R_1(z)+R_2(z))}{(1-R_1(z)-T(z))^2} = \\
&= \frac{R_1(z)^2(\RR_{11}(z) + \RR_{12}(z)+R_1(z)+R_2(z))}{z}, \\
\RR_{22}(z) &= \frac{z\RR(z)}{(1-T(z))^2} +
\frac{2z(R_1(z)+R_2(z))^2}{(1-T(z))^3} \\
&= \frac{T(z)^2\RR(z)}{z} +
\frac{2R_2(z)^2(1-T(z))}{z}.
\end{split}
\end{equation*}
Note that the system can even be solved explicitly by successive solution of
quadratic equations. We apply the method that was also used for the covariance
of $\sigma$- and $\rho$- resp. $Z$- and $\rho$-index: $\RR_{12}$ and $\RR_{22}$
are expressed in terms of the other functions. We obtain
$$\RR(z) = \RR_{11}(z) + 2\RR_{12}(z) + \RR_{22}(z) = \frac{N}{(T(z)^2-z)(R_1(z)^2-z)},$$
where the numerator $N$ is given by
\begin{multline*}
N = z\RR_{11}(z)(z+R_1(z)^2)-2R_1(z)^2R_2(z)^2(1-T(z)) \\
+2zR_1(z)^2(R_1(z)+R_2(z)) + 2zR_2(z)^2(1-T(z)).
\end{multline*}
Next, we prove that the denominator doesn't vanish within the bounds of
interest: an easy calculation shows that it can only vanish at $z = 0$, $z =
\frac{1}{4}$ or $z = \frac{4}{25}$.

So it suffices to determine the singularities of $\RR_{11}(z)$. The application
of Gr\"obner bases shows that the smallest singularity of $\RR_{11}(z)$ is the
surprisingly nice value $z_0 = \frac{8}{81}$. Expansion of $\RR_{11}$ around
$z_0$ yields
$$\RR_{11}(z) \sim \frac{2\sqrt{2}}{9} - \frac{4}{3\sqrt{7}} \sqrt{1-\frac{z}{z_0}},$$
so after calculating the values of $R_1$, $R_2$ and $T$ at $z_0$, we see that
$$\RR(z) \sim \frac{2432\sqrt{2}-1632}{3087} - \frac{32\sqrt{14}}{147}
\sqrt{1-\frac{z}{z_0}},$$
giving us the asymptotics
$$\Var(\rho_n) \sim \frac{64\sqrt{14}}{147} \cdot \left( \frac{81}{32}
\right)^n.$$

\section{The Wiener index}
\label{wiener}

Now, we want to determine the asymptotic behavior of the generating function
$$\sum_{T} X(T)W(T)z^{|T|},$$
where $X$ stands for $\sigma$-, $Z$- or $\rho$-index,$W(T)$ is the Wiener
index and $T$ runs over all rooted ordered trees. Again, we have to distinguish sets containing and not containing the
root. For the Wiener index, on the other hand, we have the following recursive
relations:
\begin{equation*}
D(T) = \sum_{i=1}^k D(T_i) + |T|-1
\end{equation*}
and
\begin{equation*}
W(T) = D(T) + \sum_{i=1}^k W(T_i) + \sum_{i \neq j} \Big(D(T_i)+|T_i|
\Big)|T_j|.
\end{equation*}
Here, $D$ is the total height or internal path length, the sum of the distances
to the root. From these relations, one finds the following functional equations
for the respective generating functions:
\begin{equation*}
\begin{split}
D(z) &= \frac{zD(z)}{(1-T(z))^2} + zT'(z) - T(z), \\
W(z) &= D(z) + \frac{zW(z)}{(1-T(z))^2} +
\frac{2z^2T'(z)(D(z)+zT'(z))}{(1-T(z))^3},
\end{split}
\end{equation*}
from which $D(z)$ and $W(z)$ can be determined without difficulty. The
asymptotic behavior of the average Wiener index follows at once:
$$E(W_n) \sim \frac{\sqrt{\pi}}{4} n^{5/2}.$$
The variance of the Wiener index is given in a paper of Janson \cite{jan03}:
$$\Var(W_n) \sim \frac{16-5\pi}{80} n^5.$$
Now, we define generating functions $\DS_1,\DS_2,\WS_1,\WS_2$ for the
product of $D(T)$ resp. $W(T)$ with the number of independent vertex subsets
containing resp. not containing the root, e.g.
$$\DS_1(z) = \sum_T D(T)\sigma_1(T) z^{|T|}.$$
In an analogous manner, we define the functions $\DZ_i,\WZ_i,\DR_i,\WR_i$.
Here, we obtain linear functional equations for the generating functions, which
can be solved explicitly.

\subsection{$\sigma$- and Wiener index}

The recursive relations give us the following system of functional equations,
which can be simplified by means of the functional equations for $S_1$ and
$S_2$, especially the facts that $\frac{z}{1-S_2(z)} = S_1(z)$,
$\frac{z}{1-S_1-S_2(z)} = S_2(z)$ and $\frac{z}{(1-S_2(z))^2} = S_2(z)$ (the
latter follows after some simple algebraic manipulations):
\begin{equation*}
\begin{split}
\DS_1(z) &= \frac{z\DS_2(z)}{(1-S_2(z))^2} + zS_1'(z) - S_1(z) \\
&= \DS_2(z)S_2(z) + zS_1'(z) - S_1(z), \\
\DS_2(z) &= \frac{z(\DS_1(z)+\DS_2(z))}{(1-S_1(z)-S_2(z))^2} + zS_2'(z) - S_2(z) \\
&= \frac{S_2(z)^2(\DS_1(z)+\DS_2(z))}{z} + zS_2'(z) - S_2(z), \\
\WS_1(z) &= \DS_1(z) + \frac{z\WS_2(z)}{(1-S_2(z))^2} +
\frac{2z^2S_2'(z)(\DS_2(z)+zS_2'(z))}{(1-S_2(z))^3} \\
&= \DS_1(z) + S_2(z)\WS_2(z) + 2S_1(z)S_2(z)S_2'(z)(\DS_2(z)+zS_2'(z)), \\
\end{split}
\end{equation*}
\begin{equation*}
\begin{split}
\WS_2(z) &= \DS_2(z) + \frac{z(\WS_1(z)+\WS_2(z))}{(1-S_1(z)-S_2(z))^2} \\
&\ \ + \frac{2z(zS_1'(z) + zS_2'(z))(\DS_1(z) + \DS_2(z)+ zS_1'(z) +
  zS_2'(z))}{(1-S_1(z) - S_2(z))^3} \\
&= \DS_2(z) + \frac{S_2(z)^2(\WS_1(z)+\WS_2(z))}{z} \\
&\ \ + \frac{2S_2(z)^3(S_1'(z) + S_2'(z))(\DS_1(z) + \DS_2(z)+ zS_1'(z) +
  zS_2'(z))}{z}.
\end{split}
\end{equation*}
All these equations are obtained as in the following example:
\begin{align*}
\DS_1(z) &= \sum_T D(T)\sigma_1(T) z^{|T|} \\
&= \sum_{k \geq 0} \sum_{T_1} \ldots \sum_{T_k} \left( \sum_{i=1}^k D(T_i)
  \prod_{j=1}^k \sigma_2(T_j) \right)z^{|T_1| + \ldots + |T_k|+1} \\
&\ \ \ + \sum_T (|T|-1)\sigma_1(T) z^{|T|} \\
&= \sum_{k \geq 0} \sum_{T_1} \ldots \sum_{T_k} \left( \sum_{i=1}^k
  D(T_i)\sigma_2(T_i) \prod_{j \neq i} \sigma_2(T_j) \right)z^{|T_1| + \ldots +
  |T_k|+1} \\
&\ \ \ + zS_1'(z) - S_1(z) \\
&= z\sum_{k \geq 0} k \DS_2(z) S_2(z)^{k-1} + zS_1'(z) - S_1(z) \\
&= \frac{z\DS_2(z)}{(1-S_2(z))^2} + zS_1'(z) - S_1(z).
\end{align*}
We solve the system of linear equations for $\WS_1$ and $\WS_2$ and obtain an
expression for $\WS(z) = \WS_1(z) + \WS_2(z)$ in terms of $S_1$ and
$S_2$. Then, we replace $S_1(z)$ by $\frac{z}{1-S_2(z)}$ and $S_1'(z)$ by
$$S_1'(z) = \frac{d}{dz} \frac{z}{1-S_2(z)} = \frac{1}{1-S_2(z)} +
\frac{zS_2'(z)}{(1-S_2(z))^2} = \frac{1}{1-S_2(z)} + S_2(z)S_2'(z).$$
Finally, implicit differentiation of the equation
$S_2(z)^3-2S_2(z)^2+S_2(z)-z=0$ yields 
$$S_2'(z) = \frac{1}{3S_2(z)^2 - 4S_2(z) + 1},$$
so $\WS$ can be written in terms of $S_2$ and $z$ only. We obtain an expression
of the form
$$\WS(z) = \frac{N}{(1-3S_2(z))^2(1-S_2(z))^3(S_2(z)^2+S_2(z)^3-z)^2},$$
where $N$ is a polynomial in $S_2$ and $z$. The denominator only vanishes at
$z=0$ (which is clearly no singularity) and at $z = \frac{4}{27}$, which is the
dominating singularity of $S_2$. Therefore, we only have to consider the
expansion of $\WS$ around $\frac{4}{27}$, which is given by
$$\WS(z) \sim \frac{5}{81\left(1-\frac{27z}{4}\right)^2}.$$
This yields
$$E(W_n\sigma_n) \sim \frac{20\sqrt{\pi}}{81} n^{5/2} \left(
  \frac{27}{16} \right)^n.$$

\subsection{$Z$- and Wiener index}

All steps are analogous to the previous section. We start from the following functional
equations, which are simplified by the relations $Z_1(z) =
\frac{z^2}{(1-Z_1(z)-Z_2(z))^3}$ and $Z_2(z) =
\frac{z}{1-Z_1(z)-Z_2(z)}$:
\begin{equation*}
\begin{split}
\DZ_1(z) &= \frac{2zZ_2(z)(\DZ_1(z)+\DZ_2(z))}{(1-Z_1(z)-Z_2(z))^3} +
\frac{z\DZ_2(z)}{(1-Z_1(z)-Z_2(z))^2} + zZ_1'(z) - Z_1(z) \\
&= \frac{2Z_1(z)Z_2(z)(\DZ_1(z)+\DZ_2(z))}{z} + \frac{Z_2(z)^2\DZ_2(z)}{z} +
zZ_1'(z) - Z_1(z), \\
\DZ_2(z) &= \frac{z(\DZ_1(z)+\DZ_2(z))}{(1-Z_1(z)-Z_2(z))^2} + zZ_2'(z) - Z_2(z) \\
&= \frac{Z_2(z)^2(\DZ_1(z)+\DZ_2(z))}{z} + zZ_2'(z) - Z_2(z), \\
\end{split}
\end{equation*}
\begin{equation*}
\begin{split}
\WZ_1(z) &= \DZ_1(z) + \frac{2zZ_2(z)(\WZ_1(z)+\WZ_2(z))}{(1-Z_1(z)-Z_2(z))^3} +
\frac{z\WZ_2(z)}{(1-Z_1(z)-Z_2(z))^2} \\
&\ \ + \frac{2z}{(1-Z_1(z)-Z_2(z))^3}
\Big((\DZ_2(z)+zZ_2'(z))(zZ_1'(z)+zZ_2'(z)) \\
&\ \ \ \ \ \ + zZ_2'(z)(\DZ_1(z)+\DZ_2(z)+zZ_1'(z)+zZ_2'(z))\Big) \\
&\ \ + \frac{6zZ_2(z)(zZ_1'(z)+zZ_2'(z))(\DZ_1(z)+\DZ_2(z)+zZ_1'(z)+zZ_2'(z))}
{(1-Z_1(z)-Z_2(z))^4} 
\\ 
&= \DZ_1(z) + \frac{2Z_1(z)Z_2(z)(\WZ_1(z)+\WZ_2(z))}{z} +
\frac{Z_2(z)^2\WZ_2(z)}{z} \\
&\ \ + 2Z_1(z) \Big((\DZ_2(z)+zZ_2'(z))(Z_1'(z)+Z_2'(z)) \\
&\ \ \ \ \ \ + Z_2'(z)(\DZ_1(z)+\DZ_2(z)+zZ_1'(z)+zZ_2'(z))\Big) \\
&\ \ + \frac{6Z_1(z)Z_2(z)^2}{z}
(Z_1'(z)+Z_2'(z))(\DZ_1(z)+\DZ_2(z)+zZ_1'(z)+zZ_2'(z)), \\ 
\WZ_2(z) &= \DZ_2(z) + \frac{z(\WZ_1(z)+\WZ_2(z))}{(1-Z_1(z)-Z_2(z))^2} \\
&\ \ + \frac{2z^2(Z_1'(z) + Z_2'(z))(\DZ_1(z) + \DZ_2(z)+ zZ_1'(z) +
  zZ_2'(z))}{(1-Z_1(z) - Z_2(z))^3} \\
&= \DZ_2(z) + \frac{Z_2(z)^2(\WZ_1(z)+\WZ_2(z))}{z} \\
&\ \ + \frac{2Z_2(z)^3(Z_1'(z) + Z_2'(z))(\DZ_1(z) + \DZ_2(z)+ zZ_1'(z) +
  zZ_2'(z))}{z}.
\end{split}
\end{equation*}
We solve this linear equation for $\WZ_1$ and $\WZ_2$ and obtain an expression
for $\WZ(z) = \WZ_1(z) + \WZ_2(z)$ in terms of $Z_1$ and $Z_2$. Then, we
replace $Z_1(z)$ by $\frac{Z_2(z)^3}{z}$ and $Z_1'(z)$ by
$\frac{3Z_2(z)^2Z_2'(z)}{z} - \frac{Z_2(z)^3}{z^2}$.
Finally, implicit differentiation of the functional equation
$z^2 - zZ_2(z) + zZ_2(z)^2 + Z_2(z)^4 = 0$ yields 
$$Z_2'(z) = \frac{2z-Z_2(z)+Z_2(z)^2}{z-2zZ_2(z)-4Z_2(z)^3},$$
so $\WZ$ can be written in terms of $Z_2$ and $z$ only. We obtain an expression
of the form
$$\WZ(z) = \frac{N}{z^3(z^2-zZ_2(z)^2-3Z_2(z)^4)^2(z-2zZ_2(z)-4Z_2(z)^3)^2},$$
where $N$ is a polynomial in $Z_2$ and $z$. The denominator only vanishes at
$z=0$ (which is clearly no singularity) and at $z =
\frac{\pm13\sqrt{13}-35}{72}$, which are singularities of $Z_2$. Therefore, we
only have to consider the expansion of $\WZ$ around the dominating
singularity $z_0 = \frac{13\sqrt{13}-35}{72}$, which is given by 
$$\WZ(z) \sim \frac{91-5\sqrt{13}}{1248\left(1-\frac{z}{z_0}\right)^2}.$$
This yields
$$E(W_nZ_n) \sim \frac{(91-5\sqrt{13})\sqrt{\pi}}{312} n^{5/2} \left(
  \frac{35+13\sqrt{13}}{54} \right)^n.$$

\subsection{$\rho$- and Wiener index}

Again, all steps are almost analogous to the previous section. We start from
the following functional equations, which are simplified by the relations
$T(z) = \frac{z}{1-T(z)}$, $R_1(z) = \frac{z}{1-R_1(z)-T(z)}$ and $R_2(z) =
\frac{z(R_1(z)+R_2(z))}{(1-T(z))^2}$:
\begin{equation*}
\begin{split}
\DR_1(z) &= \frac{z(\DR_1(z)+D(z))}{(1-R_1(z)-T(z))^2} +  zR_1'(z) - R_1(z) \\
&= \frac{R_1(z)^2(\DR_1(z)+D(z))}{z} + zR_1'(z) - R_1(z), \\
\DR_2(z) &= \frac{2zD(z)(R_1(z)+R_2(z))}{(1-T(z))^3} +
\frac{z(\DR_1(z)+\DR_2(z))}{(1-T)^2} + zR_2'(z) - R_2(z) \\
&= \frac{2D(z)T(z)R_2(z)}{z} +
\frac{T(z)^2(\DR_1(z)+\DR_2(z))}{z} + zR_2'(z) - R_2(z), \\
\WR_1(z) &= \DR_1(z) + \frac{z(\WR_1(z)+W(z))}{(1-R_1(z)+T(z))^2} \\
&\ \ + \frac{2z(zR_1'(z)+zT'(z))(\DR_1(z)+D(z)+zR_1'(z)+zT'(z))
}{(1-R_1(z)-T(z))^3} \\
&= \DR_1(z) + \frac{R_1(z)^2(\WR_1(z)+W(z))}{z} \\
&\ \ + \frac{2R_1(z)^3(R_1'(z)+T'(z))(\DR_1(z)+D(z)+zR_1'(z)+zT'(z))}{z}, \\
\WR_2(z) &= \DR_2(z) + \frac{2zW(z)(R_1(z)+R_2(z))}{(1-T(z))^3} +
\frac{z(\WR_1(z)+\WR_2(z))}{(1-T(z))^2} \\
&\ \ + \frac{2z(\DR_1(z)+\DR_2(z)+zR_1'(z)+zR_2'(z))zT'(z)}{(1-T(z))^3} \\
&\ \ + \frac{2z(zR_1'(z)+zR_2'(z))(D(z)+zT'(z))}{(1-T(z))^3} \\
&\ \ + \frac{6z^2T'(z)(D(z)+zT'(z))(R_1(z)+R_2(z))}{(1-T(z))^4} \\
&= \DR_2(z) + \frac{2T(z)^3W(z)(R_1(z)+R_2(z))}{z^2} +
\frac{T(z)^2(\WR_1(z)+\WR_2(z))}{z} \\
&\ \ + \frac{2T(z)^3(\DR_1(z)+\DR_2(z)+zR_1'(z)+zR_2'(z))T'(z)}{z} \\
&\ \ + \frac{2T(z)^3(R_1'(z)+R_2'(z))(D(z)+zT'(z))}{z} \\
&\ \ + \frac{6T(z)^4T'(z)(D(z)+zT'(z))(R_1(z)+R_2(z))}{z^2}.
\end{split}
\end{equation*}
In this case, we solve the linear equation for $\WR_1$ and $\WR_2$ and insert
the exact expressions for $T$, $D$, $W$, $R_1$ and $R_2$, which can be
determined by simple quadratic equations: we have
$$R_1(z) = \frac{1}{4}\left(1+\sqrt{1-4z}-\sqrt{2(1-10z+\sqrt{1-4z})}\right)$$
and
$$R_2(z) = \frac{1-\sqrt{1-4z}}{2\sqrt{1-4z}} \cdot R_1(z).$$
Note that we always have to take the branch whose value is 0 at $z=0$.
We will use $q_1$ as an abbreviation for $\sqrt{1-4z}$ and $q_2$ as an
abbreviation for $\sqrt{2(1-10z+\sqrt{1-4z})}$. Then we also have
$$T(z) = \frac{1-q_1}{2},\ D(z) = \frac{2z^2}{q_1^2(1+q_1)},\ W(z) =
\frac{z^2}{q_1^4}.$$
We obtain an exact
expression for $\WR(z) = \WR_1(z) + \WR_2(z)$ in terms of a rational function
in $q_1$ and $q_2$. The denominator of this expression is given by
$$q_1^6(5q_1-3)(3-5q_1+q_2)^2,$$
which only vanishes at $z=0$, $z = \frac{1}{4}$ and $z =
\frac{4}{25}$. Furthermore, $q_1$ has its only singularity at $\frac{1}{4}$,
and $q_2$ has singularities exactly at $\frac{1}{4}$ and
$\frac{4}{25}$. Therefore, the singularity of $\WR$ we have to investigate is
$\frac{4}{25}$. We obtain the expansion
$$\WR(z) \sim \frac{1}{15\left(1-\frac{25z}{4}\right)^2},$$
which yields
$$E(W_n\sigma_n) \sim \frac{4\sqrt{\pi}}{15} n^{5/2} \left(
  \frac{25}{16} \right)^n.$$

\begin{table}[htbp]
\begin{center}
\begin{tabular}{|c|}
\hline
$r(\sigma_n,Z_n) \sim (-1.01706) \cdot (0.99405)^n$ \\
$r(\sigma_n,\rho_n) \sim (1.05088) \cdot (0.99023)^n$ \\
$r(Z_n,\rho_n) \sim (-1.08924) \cdot (0.97853)^n$ \\
$r(W_n,\sigma_n) \sim (-0.27891) \cdot (0.99767)^n$ \\
$r(W_n,Z_n) \sim (0.40351) \cdot (0.99637)^n$ \\
$r(W_n,\rho_n) \sim (-1.78357) \cdot (0.98209)^n$ \\
\hline
\end{tabular}
\end{center}
\caption{Asymptotic formulas for the correlation coefficients.}
\label{asform}
\end{table}

\section{Conclusion}

From the expected values and variances which were calculated in the preceding
sections, it is possible now to deduce the asymptotic correlation coefficients
for the investigated indices (Table~\ref{asform}). For their interpretations
and further discussion see \cite{wagcor}.

\end{document}